\documentclass[twoside,12pt]{article}

\usepackage{amssymb,amsmath}

\newenvironment{proof}{\begin{trivlist}\item[]{\it
Proof.}}{\hfill$\square$\end{trivlist}}
\newenvironment{proofof}[1]{\noindent{\it Proof of
#1.}}{\hfill$\square$\\\mbox{}}

\newtheorem{theorem}{Theorem}[section]
\newtheorem{corollary}[theorem]{Corollary}

\newtheorem{lemma}[theorem]{Lemma}
\newtheorem{proposition}[theorem]{Proposition}

\newtheorem{remark}[theorem]{Remark}
\begin{document}

\date{}

\newcommand{\im}{{\rm im }}
\newcommand{\rk}{{\rm rk }}
\newcommand{\Ad}{{\rm Ad }}
\newcommand{\Tr}{{\rm Tr }}
\newcommand{\gltwo}{\mathfrak{gl}(2)}
\newcommand{\sltwo}{\mathfrak{sl}(2)}

\title{On orthogonal invariants in characteristic 2}

\author{M. Domokos ${}^{a,b}$ 
\thanks{Corresponding author. 
Supported through a European Community Marie Curie Fellowship, 
held at the University of Edinburgh. 
Partially supported by OTKA No. F 32325 and T 34530.}\; 
and P. E. Frenkel ${}^c$ 
\thanks{Partially supported by OTKA grant T 042769.}
\\ 
\\ 
{\small ${}^a$ R\'enyi Institute of Mathematics, Hungarian Academy of 
Sciences,} 
\\ {\small P.O. Box 127, 1364 Budapest, Hungary,} 
{\small E-mail: domokos@renyi.hu } 
\\ 
{\small ${}^b$ (Temporary, until January 2004) 
School of Mathematics, 
University of Edinburgh,}
\\ {\small James Clerk Maxwell Building, King's Buildings,} 
{\small Mayfield Road,} 
\\ {\small Edinburgh EH9 3JZ, Scotland}  
\\ 
\\
{\small ${}^c$ Institute of Mathematics, Budapest Technical University,} 
\\ {\small P.O.B. 91,  1521 Budapest, Hungary,}  
\\ {\small E-mail: frenkelp@renyi.hu}
}
\maketitle 
\begin{abstract} 
Working over an algebraically closed base field $k$ of characteristic $2$, 
the ring of invariants $R^G$ is studied, where $G$ is the orthogonal group 
$O(n)$ or the special orthogonal group $SO(n)$, acting naturally 
on the coordinate ring $R$ of the $m$-fold direct sum 
$k^n\oplus\cdots\oplus k^n$ of the standard vector representation. 
It is proved for $O(2)$, $O(3)=SO(3)$, $SO(4)$, and $O(4)$, 
that there exists an $m$-linear invariant with $m$ arbitrarily large, 
which is not expressible as a polynomial of invariants of lower degree. 
This is in sharp contrast with the uniform description of 
the ring of invariants valid in all other characteristics, 
and supports the conjecture that the same phenomena occur for all $n$. 
For general even $n$, new $O(n)$-invariants are constructed, 
which are not expressible as polynomials of the quadratic invariants.  
In contrast with these results, it is shown that rational invariants 
have a uniform description valid in all characteristics.  
Similarly, if $m\leq n$, then $R^{O(n)}$ is generated by the obvious 
invariants. For all $n$, the algebra $R^G$ is a finitely generated module over 
the subalgebra generated by the quadratic invariants, 
and for odd $n$, the square of any $SO(n)$-invariant is a polynomial of the 
quadratic invariants. 
Finally we mention that for even $n$, an $n$-linear $SO(n)$-invariant is 
given, which distinguishes between $SO(n)$ and $O(n)$ 
(just like the determinant in all characteristics different from $2$).
\end{abstract} 



\tableofcontents
\section{Preliminaries}\label{Prelim}
\subsection{The orthogonal group}\label{The-orth}
Let $k$ stand for an algebraically closed field of characteristic
2. Recall that the polar bilinear form $\beta$ of a quadratic form $q$ on
a finite dimensional $k$-linear space is defined by
\begin{equation}\label{betadef}\beta\left(v^{(1)},v^{(2)}\right)\overset{\rm{def}}=
q\left(v^{(1)}+v^{(2)}\right)-q\left(v^{(1)}\right)-q\left(v^{(2)}\right).\end{equation}
Note that $\beta$ is an alternating bilinear form (which implies,
but is not equivalent to, symmetry in characteristic 2). The
quadratic form $q$ is said to be non-degenerate if  $\beta
(v,\cdot)=0$  and $q(v)=0$ together imply $v=0$.

Denote coordinates in $k^n$ by $x_1$,\dots, $x_\nu$, $y_1$,\dots,
$y_\nu$ if $n=2\nu$ or  by $x_1$,\dots, $x_\nu$, $y_1$,\dots,
$y_\nu$, $z$ if $n=2\nu+1$. The orthogonal group $O(n)$ is the
group of linear isomorphisms of $k^n$ that leave the standard
non-degenerate quadratic form $$q\overset{\rm {def}}=x_1y_1 +\dots
+x_\nu y_\nu \qquad (n=2\nu) $$ resp. $$q\overset{\rm
{def}}=x_1y_1 +\dots + x_\nu y_\nu +z^2\qquad (n=2\nu +1)$$
invariant. Of course they leave the polar
form \begin{equation}\label{beta}\beta\left(v^{(1)},v^{(2)}\right)
=
x_1^{(1)} y_1^{(2)}+y_1^{(1)} x_1^{(2)}+\dots +x_\nu^{(1)}
y_\nu^{(2)}+y_\nu^{(1)} x_\nu^{(2)}\end{equation} of $q$ invariant as well. Note
that up to base change, $q$ is the only non-degenerate quadratic form on $k^n$.

The form $\beta$ is non-degenerate if and only if $n$ is even. For
$n=2\nu +1$, $$\ker\beta\overset{\rm {def}}=\left\{v:\beta
(v,w)=0\;\; {\rm {for}\; \rm {all }}\; w\right\}$$ is the $z$
axis.

The symplectic group $Sp(2\nu)$ is the group of linear
isomorphisms of $k^{2\nu}$ that leave the standard symplectic form
$\beta$ invariant. So $O(2\nu)\leq Sp(2\nu)\leq SL(2\nu)$. In fact
$O(n)\leq SL(n)$ for all $n$. The algebraic group $O(n)$ is
connected for odd $n$ and has two components for even $n$.  For
all $n$, the component containing the identity is the special
orthogonal group  $SO(n)$ (this can be taken as the definition of
$SO(n)$). Thus, $SO(2\nu+1)=O(2\nu+1)$, whereas $SO(2\nu)$ is a subgroup of
index 2 in $O(2\nu)$. 

Call a vector $u$ {\it non-singular} if $q(u)\neq 0$. 
For a non-singular vector $u$, we write $T_u$ for the {\it reflection} 
defined by 
$$T_uv\overset{\rm def}=v-\frac{\beta
(v,u)}{q(u)}u.$$
It is well known that $O(n)$ is generated by
reflections, and $SO(n)$ is the set of elements that are
expressible as a product of an even number of reflections.  

For $n=2\nu +1$, each $A\in O(2\nu +1)$ acts as the identity on
the $z$ axis and acts symplectically on the factor space $k^{2\nu
+1}/\ker\beta$. This gives a homomorphism $\phi:O(2\nu +1)\to
Sp(2\nu)$ which is in fact an isomorphism (of groups, but not of
algebraic groups). See \cite[Theorem 11.9]{T} for a proof.

\subsection{Invariants}
We write $R$ or $R_{n\times m}$ for the algebra of polynomials in
the coordinates of the indeterminate $n$-dimensional vectors
$v^{(1)}$, \dots, $v^{(m)}$. We write $K$ or $K_{n\times m}$ for the
field of rational functions. A $G$ in the superscript indicates
the subalgebra (sub-field) formed by the functions invariant under
the subgroup $G$ of $GL(n)$ acting on $m$-tuples of vectors in the
obvious way. Let
\begin{equation}\label{alapinv}
\begin{aligned}
Q^{(i)}&\overset{\rm{def}}=q\left(v^{(i)}\right) & =&
\begin{cases}x_1^{(i)}y_1^{(i)}+ \dots +x_\nu^{(i)}y_\nu^{(i)}\\
x_1^{(i)}y_1^{(i)}+ \dots
+x_\nu^{(i)}y_\nu^{(i)}+{z^{(i)}}^2,\end{cases}\\
B^{(ij)}&\overset{\rm{def}}=\beta\left(v^{(i)},v^{(j)}\right) & =&
x_1^{(i)}y_1^{(j)}+y_1^{(i)}x_1^{(j)}+ \dots
+x_\nu^{(i)}y_\nu^{(j)}+y_\nu^{(i)}x_\nu^{(j)}.\end{aligned}
\end{equation} Let $$
D^{(i_1, \dots , i_n)}\overset{\rm{def}}=\det \left[v^{(i_1)},
\dots , v^{(i_n)}\right]$$ be the determinant of the matrix that
has $v^{(i_1)}$, \dots, $v^{(i_n)}$ as its columns. Then  $
Q^{(i)}$, $B^{(ij)}$, $D^{(i_1, \dots , i_n)}$ are
multi-homogeneous elements of $ R_{n\times m}^{O(n)}$.

(By the multi-degree of a monomial in the polynomial ring $R_{n\times m}$ 
we mean $\alpha=(\alpha^{(1)},\dots,\alpha^{(m)})$, where $\alpha^{(i)}$ 
is the total degree of the monomial in the variables belonging to $v^{(i)}$. 
The action of $O(n)$ preserves this multi-degree, therefore, 
$R_{n\times m}^{O(n)}$ is spanned by multi-homogeneous elements. 
A multi-homogeneous invariant of multi-degree $(1,\dots,1)$ will be called 
{\it multi-linear}.)

It is a classical fact that over a field of characteristic zero,
the algebra $R^{O(n)}$ is generated by the scalar products $B^{(ij)}$ of the
indeterminate vectors under consideration, and the algebra
$R^{SO(n)}$ is generated by the scalar products and the
determinants. That is  the so-called ``first fundamental theorem''
for the (special) orthogonal group; it has been discussed along
with the analogous results for the other  classical groups in
Hermann Weyl's work \cite{W}. De Concini and Procesi \cite{CP} gave a characteristic free treatment 
to the subject, in particular, they proved that the first fundamental 
theorem for the (special) orthogonal group remains unchanged 
in odd characteristic. 
Concerning characteristic $2$, Richman \cite{Ri} proved later 
that the algebra $R^G$ for the group $G$ 
preserving the bilinear form 
$x_1^{(1)}x_1^{(2)}+\cdots+x_n^{(1)}x_n^{(2)}$ 
is generated in degree 
$1$ and $2$. 
However, though this group preserves the quadratic form 
$x_1^2+\cdots+x_n^2$, it is not the so-called `orthogonal group' 
in characteristic $2$: the quadratic form $x_1^2+\cdots+x_n^2$ 
is the square of a linear form, hence is degenerate. So the question 
about vector invariants of the orthogonal group remains open 
in characteristic $2$, when the behaviour of invariants 
turns out to be very much different, see Section~\ref{high}.

\section{Indecomposable invariants of high degree}\label{high}
The results in this section  make the following conjecture
plausible: for any fixed $n\geq 2$ (resp. $n\geq 3$), there exist
arbitrarily large values of $m$ and $m$-linear invariants $f_m\in
R_{n\times m}^{O(n)}$ (resp. $R_{n\times m}^{SO(n)}$) such that
$f_m$ cannot be expressed as a polynomial in invariants of lower
degree.
We  prove this  for $n\leq 4$. The paper \cite{DKZ} contained a
more sophisticated proof for the $SO(4)$ case. It was first
pointed out in this paper that special orthogonal invariants
behave much differently in characteristic 2.

In the general case, we have no proof of the conjecture, but in
Subsection~\ref{even} we shall prove at least that the algebra
$R_{n\times m}^{O(n)}$ is not generated by the $Q^{(i)}$ and
$B^{(ij)}$ if $n\geq 2$ and $m$ is large enough (compared to $n$).
This is obvious for odd $n$, since if $m\geq n$, then $D^{(1\cdots
n)}$ is not expressible as a polynomial in the $Q^{(i)}$ and
$B^{(ij)}$, but it is non-trivial for even $n$.

\subsection{The two-dimensional case}
To treat the two-dimensional case, observe that the matrix $$A=
 \left(
\begin{matrix}
a_{11}  &  a_{12} \\ a_{21}  &  a_{22}
\end{matrix}
\right)$$ is orthogonal if and only if
$$(a_{11}x+a_{12}y)(a_{21}x+a_{22}y)=xy,$$ that is,
$a_{11}a_{21}=a_{12}a_{22}=0$ and $a_{11}a_{22}+a_{12}a_{21}=1$.
So $$ O(2)=\left\{ \left(
\begin{matrix}
a  & 0 \\ 0 & 1/a \end{matrix} \right): a\in k^* \right\}\cup
\left\{ \left(
\begin{matrix}
0  & a \\ 1/a & 0 \end{matrix} \right): a\in k^* \right\},$$
 where the first of the two terms is $SO(2)$.

Therefore, a polynomial is invariant under $SO(2)$ if and only if
all its terms have the same number of   $x$'s and $y$'s. It
follows that the algebra of $SO(2)$-invariant polynomials is
generated by quadratic elements: $$ R_{2\times
m}^{SO(2)}=k\left[x^{(i)}y^{(j)}:i,j= 1, \dots , m\right].$$

That is not the case with $O(2)$-invariants. An $SO(2)$-invariant
 is $O(2)$-invariant exactly if it is invariant under $
\left(
\begin {matrix}
0 & 1 \\ 1 & 0
\end{matrix}
\right)$, that is, exactly if it is a linear combination over $k$
of (multi-homogen\-eous) polynomials of the form $$
x^{(i_1)}\cdots x^{(i_s)}y^{(i_1)}\cdots y^{(i_s)}=Q^{(i_1)}\cdots
Q^{(i_s)}$$ and $$x^{(i_1)}\cdots x^{(i_s)}y^{(j_1)}\cdots
y^{(j_s)}+y^{(i_1)}\cdots y^{(i_s)}x^{(j_1)}\cdots x^{(j_s)}
\overset{\text {def}} =B^{(i_1, \dots, i_s | j_1, \dots, j_s)}.$$
(Note that the new notation is in accordance with the notation
$B^{(ij)}$ introduced before.)

\begin{proposition}\label{2dim}
\begin{itemize}
\item[(i)] Assume that the indices  $i_1$, \dots ,  $i_s$, $j_1$, \dots ,
$j_s$  are all different. Then the $O(2)$-invariant $B^{(i_1,
\dots, i_s|j_1, \dots, j_s)}$ is not expressible as a polynomial
in invariants of lower degree.

\item[(ii)] Assume that each of the indices  $1$, \dots,  $m$  occurs
among the indices  $i_1$, \dots, $i_s$, $j_1$, \dots, $j_s$  the
same number of times as it occurs among the indices  $i_1'$,
\dots, $i_s'$, $j_1'$, \dots, $j_s'$. Then the multi-homogen\-eous
$O(2)$-invariant $B^{(i_1, \dots, i_s|j_1, \dots, j_s)}+B^{(i_1',
\dots, i_s'|j_1', \dots, j_s')}$ (if non-zero) is expressible as
the product of two $B$'s of lower degree.

\item[(iii)]  Assume that the indices  $i_1$, \dots,  $i_s$, $j_1$,
\dots, $j_s$  are not all different. Then the $O(2)$-invariant
$B^{(i_1, \dots, i_s|j_1, \dots, j_s)}$ is  expressible as a
polynomial in invariants of lower degree.
\end{itemize}
\end{proposition}

\begin{proof}
(i) Let $ \alpha$ denote the multi-degree of $B^{(i_1, \dots,
i_s|j_1, \dots, j_s)}$. So $\alpha^{(i)}=1$ if $i$ is one of the
indices $i_1$, \dots , $i_s$, $j_1$, \dots, $j_s$; and
$\alpha^{(i)}=0$ otherwise.

We only need to prove that $B^{(i_1, \dots, i_s|j_1, \dots, j_s)}$
is not expressible as a linear combination of products with two
factors each, each factor being of lower multi-degree and being
either some $B$
or some product of $Q$'s. Note that such a  product (of two
factors) is always multi-homogeneous; its multi-degree is $\mathbb
\alpha$ if and only if both factors are $B$'s (of lower
multi-degree) with no repetition of indices and with $\{i_1, \dots
,i_s,j_1, \dots ,j_s\}$ as the disjoint union of the two
index-sets. But such a product is always the sum of \it two \rm
$B$'s of multi-degree $\mathbb \alpha$. Therefore, any linear
combination of such products, when expressed as a linear
combination of the $B$'s of multi-degree $\mathbb \alpha$, gives
rise to coefficients that add up to zero.
The statement follows.

\smallskip

(ii) The assumption can be formulated by writing $I+J=I'+J'$ for
the multi-sets $I=\{i_1, \dots , i_s\}$, $J=\{j_1, \dots , j_s\}$,
$I'=\{i_1', \dots , i_s'\}$, and $J'=\{j_1', \dots , j_s'\}$. It
follows that $I=E+G$, $J=F+H$, $I'=E+H$, and $J'=F+G$ with
suitable multi-sets $E$, $F$, $G$, and $H$. That implies
$|E|=|F|$, $|G|=|H|$, and $$
B^{(E|F)}B^{(G|H)}=B^{(E+G|F+H)}+B^{(E+H|F+G)}=B^{(I|J)}+B^{(I'|J')}.$$

\smallskip

(iii) Using (ii), we may assume  $i_1=j_1$. Then $$B^{(i_1, \dots
,i_s|j_1, \dots, j_s)}=Q^{(i_1)}B^{(i_2, \dots ,i_s|j_2, \dots
,j_s)}.$$
\end {proof}

The following theorem is an easy consequence.

\begin{theorem}\label{fund}
\begin{itemize}
\item[(i)] The algebra $R_{2\times m}^{O(2)}$ is generated by the
invariants $$Q^{(i)} \qquad\textrm{ and }\qquad B^{(i_1, \dots ,i_s|j_1, \dots ,j_s)},$$
where $1\leq i \leq m$ and $1\leq i_1< \dots <i_s<j_1< \dots
<j_s\leq m$, respectively.

\item[(ii)] The system of generators in (i) is minimal. Indeed, any
system of multi-homogen\-eous generators of the algebra
$R_{2\times m}^{O(2)}$ must contain the invariants $Q^{(i)}$
(possibly multiplied by non-zero constants), and must contain
invariants of multi-degree  $\mathbb \alpha$ for all 0-1 sequences
$\mathbb \alpha=\left(\alpha^{(1)}, \dots , \alpha^{(m)}\right)$
that contain an even number of 1's.
\end{itemize}
\end{theorem}

\subsection{The three-dimensional case}\label{three}
Let us interpret  $k^3$ as $\mathfrak {\mathfrak{sl}}(2)$ via
$$v=\left(\begin{matrix}x\\y\\z\end{matrix}\right) \leftrightarrow
\left(
\begin{matrix}  z  & x \\ y  & z
\end{matrix}\right)=V.$$ Then
$q(v)=xy+z^2=\det V$. So, for any $T\in SL(2)$, $$ \Ad\;
T:\mathfrak{sl}(2) \to \mathfrak {sl}(2), \qquad V\mapsto
TVT^{-1}$$ is orthogonal. It is easily seen that every orthogonal
transformation is of this form. So, for $i_1, \dots , i_s\in \{1,
\dots ,  m\}$, the polynomial

$${\textit Tr}^{(i_1, \dots, i_s)} \overset {\rm def}=\Tr
\left(V^{(i_1)}\cdots V^{(i_s)}\right)$$ is $O(3)$-invariant:
 ${\textit Tr}^{(i_1, \dots, i_s)} \in R_{3\times m}^{O(3)}$.

\begin{proposition}\label{Tr}
If the indices $i_1$, \dots , $i_s$ are all different, then $
{\textit Tr}^{(i_1, \dots, i_s)}$ is not expressible as a
polynomial in $O(3)$-invariants of lower degree.
\end{proposition}

\begin{proof}
 We may assume $s=m$ and $i_1=1$, \dots ,
$i_s=s$. We first assume that $s$ is even; say, $s=2\sigma$.

If we replace every occurrence of all  the variables $z^{(i)}$ in
an $O(3)$-invariant by zero, then we get an  $O(2)$-invariant,
since $A\oplus 1\in O(3)$
 if $A\in O(2)$. The degree is unchanged or decreased. Therefore, it suffices to prove that
the $O(2)$-invariant
\begin{align*} {\textit Tr}^{(1, \dots, 2\sigma)}\mid_{z=0}
=\Tr \left(  \left(\begin{matrix} 0 & x^{(1)}\\y^{(1)} &
0\end{matrix} \right)
  \left(\begin{matrix}0 & x^{(2)}\\y^{(2)} & 0\end{matrix}\right)
  \cdots
\left(\begin{matrix}0 & x^{(2\sigma-1)}\\y^{(2\sigma-1)} &
0\end{matrix}\right)
  \left(\begin{matrix}0 & x^{(2\sigma)}\\y^{(2\sigma)} & 0\end{matrix}\right)
\right)
=
\\
 =\Tr\left(
  \left(
\begin{matrix}x^{(1)}y^{(2)} & \\ & y^{(1)}x^{(2)}\end{matrix}\right)
 \cdots  {\left(
\begin{matrix}x^{(2\sigma-1)}y^{(2\sigma)} & \\ & y^{(2\sigma-1)}x^{(2\sigma)}\end{matrix}\right) } \right)
=B^{(1,3, \dots, 2\sigma-1|2,4, \dots, 2\sigma)}
\end{align*}
is not expressible as a polynomial in $O(2)$-invariants of lower
degree. That was the statement of  Proposition~\ref{2dim}(i).

Assume now that $s$ is odd. Assume indirectly that
\begin{equation}
\label{eq:decomposition} Tr^{(1,\ldots,s)}=\sum_j a_jb_j,
\end{equation}
where $a_j,b_j\in R^{O(3)}$ are multi-homogeneous invariants of
strictly positive degree. We may assume that $a_jb_j$ is
$s$-linear for all $j$. Denote by $\pi:R_{3\times s}\to R_{3\times
(s-1)}$ the algebra homomorphism induced by the embedding
$\sltwo^{m-1}\to \sltwo^m$,
$\left(V^{(1)},\ldots,V^{(s-1)}\right)\mapsto
\left(V^{(1)},\ldots,V^{(s-1)},I\right)$, where $I$ stands for the
identity matrix. 
Obviously, $\pi$ maps a multi-homogeneous polynomial of 
multi-degree $(\alpha_1,\ldots,\alpha_{s-1},\alpha_s)$ 
to a multi-homogeneous polynomial of multi-degree 
$(\alpha_1,\ldots,\alpha_{s-1})$. 
Since $I$ is fixed by the $SL(2)$-action, we have
that $\pi(R_{3\times s}^{O(3)})\subseteq R_{3\times
{(s-1)}}^{O(3)}$. Applying the map $\pi$ to
\eqref{eq:decomposition} we get that
$Tr^{(1,\ldots,s-1)}=\sum\pi(a_j)\pi(b_j)$. The degree of $a_j$
and $b_j$ is at least $2$ for all $j$, hence each  $\pi(a_j)$ and
each $\pi(b_j)$ is either zero or a homogeneous $O(3)$-invariant
of positive degree. Therefore, $Tr^{(1,\ldots,s-1)}$ can be
expressed by  invariants of lower degree. But $s-1$ is even, so
this contradicts what we have already proven.
\end{proof}

\begin{theorem} \label{2.4}
A minimal system of generators of $R_{3\times m}^{O(3)}$ is
\[\{Q^{(j)},\ {\textit Tr}^{(i_1,\ldots,i_s)}\mid
1\leq j\leq m; \ 2\leq s\leq m; \ 1\leq i_1<\cdots <i_s\leq m\}.\]
\end{theorem}

\begin{proof}
The group $SL(2)$ acts on $\gltwo$, the space of $2\times 2$
matrices, by conjugation. Denote by $\gltwo^m$ (resp. $\sltwo^m$)
the $m$-fold direct sum of copies of $\gltwo$ (resp. $\sltwo$),
endowed with the diagonal $SL(2)$-action. Denote by $P$ the
coordinate ring of $\gltwo^m$, and recall that $R$ is the
coordinate ring of $\sltwo^m$. Restriction of functions from
$\gltwo^m$ to $\sltwo^m$ induces a surjective algebra homomorphism
$\varphi:P\to R$. Clearly we have $\varphi(P^{SL(2)})\subseteq
R^{O(3)}$. Now $(\gltwo^m,\sltwo^m)$ is a good pair of
$SL(2)$-varieties in the sense of \cite{D1990}. This follows for
example from \cite[Proposition 1.3b]{D1990}, since $\sltwo^m$ is
an $m$-codimensional linear subspace in the good variety
$\gltwo^m$, defined as the zero locus of $m$ linear
$SL(2)$-invariants on $\gltwo^m$, hence $\sltwo^m$ is a good
complete intersection in $\gltwo^m$. As a consequence of general
properties of modules with good filtrations (cf. \cite{D1990}) we
get that the restriction of $\varphi$ to the ring of invariants
$P^{SL(2)}$ is surjective onto $R^{O(3)}$. In particular, a
generating system of $P^{SL(2)}$ is mapped to a generating system
of $R^{O(3)}$. Using the result of \cite{D1992}, a minimal system
of generators of $P^{SL(2)}$ was determined in \cite{DKZ}. This is
mapped by $\varphi$ to the generating system of $R^{O(3)}$ stated
in our theorem. So the only thing left to show is that the above
generating system is minimal. Since it consists of
multi-homogeneous elements with pairwise different multi-degree,
it is sufficient to prove that none of them can be expressed by
invariants of strictly lower degree. This is clear for $Q^{(j)}$,
and this is the content of Proposition~\ref{Tr} for
$Tr^{(i_1,\ldots,i_s)}$.
\end{proof}

\subsection{The four-dimensional case}\label{four}
To treat the four-dimensional case, we interpret $k^4$ as
$\mathfrak {\mathfrak{gl}}(2)$ via $$v=\left(\begin{matrix}
x_1\\x_2\\y_1\\y_2\end{matrix}\right) \leftrightarrow\left(
\begin{matrix}  x_1  & x_2 \\ y_2  & y_1
\end{matrix}\right)=V.$$ Then
$q(v)=x_1y_1+x_2y_2=\det V$. So, for any $S,T\in SL(2)$, the
transformation $$ \mathfrak{gl}(2) \to \mathfrak {gl}(2), \qquad
V\mapsto SVT^{-1}$$ is orthogonal.  We get a homomorphism
$\varphi:SL(2)\times SL(2)\to O(4)$ which is easily seen to be
injective; its image is a six-dimensional irreducible subgroup of
$O(4)$, so it must be $SO(4)$. This interpretation of $SO(4)$ was
used in \cite{DKZ} to show that the algebra $R_{4\times
m}^{SO(4)}$ is not generated by its elements of degree $<m-1$.
A simpler proof can be given by means of the
following construction.

Let $i_1, \dots, i_s;j_1, \dots, j_s\in \{1 \dots m\}$. The
determinant $$\left|
\begin{matrix}
V^{(i_1)} & V^{(j_1)} &  & &
\\ & V^{(i_2)} & V^{(j_2)} &  & \\ &  & \ddots &\ddots  & \\ &  &
 & V^{(i_{s-1})} & V^{(j_{s-1})} \\ V^{(j_s)} &  &  &  & V^{(i_s)}
\end{matrix}
\right|$$ is $SO(4)$-invariant. Assume that the indices $i_1$,
\dots , $i_s$; $j_1$ \dots,  $j_s$ are all different. The
$2s$-linear component of the above determinant is also invariant,
denote it by $$F=F^{(i_1, \dots, i_{s}|j_1 \dots, j_s)}\in
R_{4\times m}^{SO(4)}.$$

\begin{proposition}\label{F}
If the indices $i_1$, \dots ,  $i_{s}$; $j_1$, \dots,
$j_s\in\{1,\dots,m\}$ are all different, then $ F ^{(i_1, \dots,
i_{s}|j_1, \dots, j_s)}$ is not expressible as a polynomial in
$SO(4)$-invariants of lower degree.
\end{proposition}

\begin{proof}
If we replace every occurrence of all  the variables $x_2^{(i)}$
and $y_2^{(i)}$ in an $SO(4)$-invariant by zero, then we get an
$O(2)$-invariant, since  if $A\in O(2)$ then $A\oplus
\left(\begin{matrix}1&0\\0&1\end{matrix}\right)$ and
$A\oplus\left(\begin{matrix}0&1\\1&0\end{matrix}\right)$ are both
in $O(4)$ 
(where $A$ acts on the $x_1,y_1$ coordinate plane) 
and one of them must be in $SO(4)$. The degree is
unchanged or decreased. Therefore, it suffices to prove that the
$O(2)$-invariant
 $$F^{(i_1, \dots,
i_s|j_1, \dots, j_s)}\mid_{x_2=y_2=0}$$ is not expressible as a
polynomial in  $O(2)$-invariants of degree $<2s$. That is the
$2s$-linear component of the determinant $$\left|\begin{matrix}
x^{(i_1)}&&x^{(j_1)}&&&&&&&\\ &y^{(i_1)}&&y^{(j_1)}&&&&&&\\
&&x^{(i_2)}&&x^{(j_2)}&&&&&\\ &&&y^{(i_2)}&&y^{(j_2)}&&&&\\
&&&&\ddots&&\ddots&&&\\ &&&&&\ddots&&\ddots&&\\
&&&&&&x^{(i_{s-1})}&&x^{(j_{s-1})}&\\
&&&&&&&y^{(i_{s-1})}&&y^{(j_{s-1})}\\
x^{(j_s)}&&&&&&&&x^{(i_s)}&\\
&y^{(j_s)}&&&&&&&&y^{(i_s)}\end{matrix}\right|,$$
which is nothing but $B^{(i_1, \dots, i_s|j_1, \dots, j_s)}$. 
The statement follows from Proposition~\ref{2dim}(i).
\end{proof}

\begin{corollary}
Any system of multi-homogen\-eous generators of the algebra
$R_{4\times m}^{SO(4)}$ must contain the invariants $Q^{(i)}$
(possibly multiplied by non-zero constants), and must contain
invariants of multi-degree  $\mathbb \alpha$ for all 0-1 sequences
$\mathbb \alpha=\left(\alpha^{(1)}, \dots , \alpha^{(m)}\right)$
that contain an even number of 1's.
\end{corollary}

To treat the full orthogonal group $O(4)$, consider the sum
$$G=G^{(i_1, \dots, i_s|j_1, \dots, j_s)}=F^{(i_1, \dots, i_s|j_1,
\dots, j_s)}+\sigma F^{(i_1, \dots, i_s|j_1, \dots, j_s)},$$ where
$\sigma$ represents the coset $O(4)\backslash SO(4)$. Obviously,
$G$ is invariant under $O(4)$.

\begin{proposition}\label{G}
If the indices $i_1$, \dots ,  $i_{s}$; $j_1$, \dots,
$j_s\in\{1,\dots,m\}$ are all different, then $ G ^{(i_1, \dots,
i_{s}|j_1, \dots, j_s)}$ is not expressible as a polynomial in
$O(4)$-invariants of  degree less than $2s-2$.
\end{proposition}

\begin{proof}
The substitution \begin{align*}x_2^{(i_t)}=y_2^{(i_t)}=0\qquad
&(t=1,\dots,s),\\x_2^{(j_t)}=y_2^{(j_t)}=0\qquad
&(t=1,\dots,s-2),\\V^{(j_{s-1})}=\left(\begin{matrix} 0 &
0\\1&0\end{matrix}\right),\qquad &V^{(j_{s})}=\left(\begin{matrix}
0 & 1\\0&0\end{matrix}\right)\end{align*} turns any
$O(4)$-invariant into an $O(2)$-invariant, since 
we may embed $O(2)$ into $O(4)$ by identifying $A\in O(2)$
with $A\oplus \left(\begin{matrix}1&0\\0&1\end{matrix}\right)\in
O(4)$ (where $A$ acts on the $x_1,y_1$ coordinate plane), 
and the subspace of $k^{4\times m}$ defined by the above equations is 
stable under $O(2)$. 
The degree is unchanged or decreased. Therefore, it
suffices to prove that the $O(2)$-invariant that $ G $ 
turns into is not expressible as a
polynomial in $O(2)$-invariants of degree $<2s-2$. Now $ F $ 
turns into  the $(2s-2)$-linear component of the determinant
$$\left|\begin{matrix} x^{(i_1)}&
&x^{(j_1)}&
&&&&&&\\
&y^{(i_1)}&
&y^{(j_1)}&&&&&&\\ &&x^{(i_2)}&
&x^{(j_2)}&
&&&&\\
&&
&y^{(i_2)}&
&y^{(j_2)}&&&&\\ &&&&\ddots&&\ddots&&&\\
&&&&&\ddots&&\ddots&&\\ &&&&&&x^{(i_{s-1})}&
&0 &0
\\
&&&&&&
&y^{(i_{s-1})}&1&0
\\ 0
&1&&&&&&&x^{(i_s)}&
\\0
&0
&&&&&&&
&y^{(i_s)}\end{matrix}\right|,$$
which is nothing but $x^{(i_1)}x^{(i_2)}\cdots
x^{(i_{s-1})}y^{(i_s)}y^{(j_1)}\cdots y^{(j_{s-2})}$. Since the
representative $\sigma: x_1\leftrightarrow y_1$ of $O(4)\backslash
SO(4)$ commutes with the substitution under consideration, $\sigma
F$ is turned  into $y^{(i_1)}y^{(i_2)}\cdots
y^{(i_{s-1})}x^{(i_s)}x^{(j_1)}\cdots x^{(j_{s-2})}$ and therefore
$G$ is turned into $B^{(i_1,i_2, \dots, i_{s-1} | i_s,j_1, \dots,
j_{s-2})}$.
The statement follows from Proposition~\ref{2dim}(i).
\end{proof}

\subsection{The even-dimensional case}\label{even}
We turn to the even-dimensional case in general. To a monomial
depending on the vectors $v^{(1)}$, \dots, $v^{(m)}$ in a multi-linear
fashion we shall associate the $2\times\nu$ matrix
$$\left(\begin{matrix}\sigma_1 & \cdots & \sigma_\nu\\\tau_1 &
\cdots & \tau_\nu\end{matrix}\right)$$ called the \it type \rm of
the monomial, whose entry $\sigma_t$ is the number of occurrences
of $x_t$ as a factor of the monomial, and $\tau_t$ is the number
of occurrences of $y_t$. So $$m=\sigma_1+\tau_1+ \dots
+\sigma_\nu+\tau_\nu.$$

\begin{lemma}\label{pinv}
Denote by $p$ the sum of all  monomials that depend on the
two-di\-men\-sion\-al vectors $v^{(1)}$, \dots, $v^{(6)}$ in a
sextilinear fashion and have type
$\left(\begin{matrix}3\\3\end{matrix}\right)$. Then $p$ is a
unimodular invariant: $p\in R_{2\times  6}^{SL(2)}$.
\end{lemma}

\begin{proof}
Invariance under $$\left(\begin{matrix} c & \\ &1/c
\end{matrix}\right)
\qquad \left(c\in k^*\right) $$ is obvious as all terms in $p$ are
invariant. It suffices to check invariance under $$A_c= \left(
\begin{matrix}
1 & c\\  &1
\end{matrix}
\right) \qquad{\textrm {and}}\qquad A_c^T= \left(
\begin{matrix}
1 & \\ c&1
\end{matrix}
\right) \qquad\qquad (c \in k).$$ By symmetry, it is sufficient to
deal with  $A_c$. The transformed polynomial $$ p_c\left(v^{(1)},
\dots,  v^{(6)}\right)= p\left(A_cv^{(1)}, \dots,
A_cv^{(6)}\right) $$ is a linear combination of sextilinear
monomials whose type
$\left(\begin{matrix}\sigma\\\tau\end{matrix}\right)$ (where
$\sigma+\tau=6$) satisfies the inequality $\sigma\leq 3\leq\tau$.
The coefficient of such a monomial is $c^{\tau-3}\binom{\tau}{3}$.
That is 1 if $\tau=3$ and zero otherwise. So $p_c=p$.
\end{proof}

Let $m\geq 2(3\nu -1)$, and denote by  $f\in R_{2\nu\times 2(3\nu
-1)}\leq R_{2\nu\times m}$ the sum of all monomials that depend in
a multi-linear fashion on the first $2(3\nu -1)$ indeterminate
vectors (and do not involve the rest), and the two rows of whose
type coincide, each row being a permutation of $(2,3,3,\dots,3)$ (one 2 and $(\nu -1)$ 3's).

\begin{theorem}
\begin{itemize}
\item[(i)] The polynomial $f$ is an orthogonal invariant.

\item[(ii)] The polynomial $f$ is not expressible as a polynomial in the
$Q^{(i)}$ and $B^{(ij)}$.
\end{itemize}
\end{theorem}

\begin{proof}
(i) Let $A^{-T}$ denote the inverse transpose  of the matrix $A$.
It suffices to check invariance under the subgroup formed by
transformations of the form $$\hat A=\left(\begin{matrix}A &
\\ & A^{-T}\end{matrix}\right) \qquad (A\in GL(\nu))$$ and under
the reflection $x_1\leftrightarrow y_1$, as these generate
$O(2\nu)$. (This follows easily from the fact that $O(2\nu)$ is
generated by reflections. Indeed, $x_1\leftrightarrow y_1$ can be
turned into an  arbitrary reflection via conjugation by some $\hat
A$, since $\{\hat{A}\mid A\in GL(\nu)\}$ acts transitively on 
$\{v\in k^n\mid q(v)=1\}$.)

Invariance under $x_1\leftrightarrow y_1$ is obvious as the terms
in $f$ simply undergo a permutation (of order 2). Now look at
$\hat A$. We may restrict $A$ to a system of generators of
$GL(\nu)$.

Invariance under $$\hat A_{i,c}: x_i\mapsto cx_i, \quad y_i\mapsto
c^{-1}y_i\qquad (i\in\{1, \dots, m\},\; c\in k^*)$$ is obvious as
each term in $f$ is invariant.

By symmetry, it suffices to check invariance under  $$\hat A_c:
x_1\mapsto x_1+cx_2,\quad y_2\mapsto cy_1+y_2\qquad (c\in k). $$
To this end, write $f$ as $f=g+h$ where $g$ is the sum of those
terms in $f$ that have only 3's in the first two columns of their
types, and $h$ is the sum of the other terms.

We use the above lemma to prove that $g$ is in fact invariant not
just under $\hat A_c$, but under both of the transformations $
x_1\mapsto x_1+cx_2$ and $y_2\mapsto cy_1+y_2$. By symmetry, it is
sufficient to deal with $ x_1\mapsto x_1+cx_2$. Let us break  $g$
up into sub-sums in the following way. Two terms shall be in the
same sub-sum if and only if the six vector variables whose $x_1$
or $x_2$ coordinate is involved are the same for the two terms,
and each of the other $2(3\nu -1)-6$ vector variables involved is
involved in the two terms via the same coordinate. Each sub-sum
will then consist of $\binom{6}{3}$ terms whose sum is invariant
under $ x_1\mapsto x_1+cx_2$ by the above lemma.

We are left with the task of proving that $h$ is invariant under
$\hat A_c$. To this end, let us break $h$ up into sub-sums in the
following way. Two terms shall be in the same sub-sum if and only
if the ten vector variables whose   $x_1$, $x_2$, $y_1$ or $y_2$
coordinate is involved are the same for the two terms, and each of
the other $2(3\nu -1)-10$ vector variables involved is involved in
the two terms via the same coordinate. Each sub-sum will then
consist of $2\binom{10}{5}\binom{5}{2}^2$ terms whose sum is
invariant under $\hat A_c$, as we shall now check. In other words,
we have to check that the sum $r$ of all monomials that depend in
a decilinear fashion on the four-dimensional vectors $v^{(1)}$,
\dots, $v^{(10)}$ and have type
$\left(\begin{matrix}2&3\\2&3\end{matrix}\right)$ or
$\left(\begin{matrix}3&2\\3&2\end{matrix}\right)$ is  invariant
under $\hat A_c$. The transformed polynomial $$ r_c\left(v^{(1)},
\dots,  v^{(10)}\right)= r\left(\hat A_cv^{(1)}, \dots,  \hat
A_cv^{(10)}\right) $$ is a linear combination of decilinear
monomials whose type $\left(\begin{matrix}\sigma_1 &
\sigma_2\\\tau_1 & \tau_2
\end{matrix}\right)$ satisfies
$\sigma_1+\sigma_2=\tau_1+\tau_2=5$. The coefficient of such a
monomial is
$$
c^{\sigma_2-3+\tau_1-2}\binom{\sigma_2}{3}\binom{\tau_1}{2}+c^{\tau_1-3+\sigma_2-2}\binom{\tau_1}{3}\binom{\sigma_2}{2}
=c^{\sigma_2+\tau_1-5}\left(\binom{\sigma_2}{3}\binom{\tau_1}{2}+\binom{\tau_1}{3}\binom{\sigma_2}{2}\right).$$
That is 1 if $\{\tau_1,\sigma_2\}=\{2,3\}$ and zero otherwise. So
$r_c=r$.

\smallskip

(ii) Since $f$ is multi-linear, the only way for the proposition
to be false would be if $f$ were a polynomial in the $B^{(ij)}$. So
it suffices to show that $f$ is not a symplectic invariant. We
show that it is not invariant under the symplectic transformation
$$T:x_1\mapsto x_1+y_1.$$  Write $f$ as  $f=\tilde g+\tilde h$
where $\tilde g$ is the sum of those terms in $f$ that have three
$x_1$'s and three  $y_1$'s among their factors, and $\tilde h$ is
the sum of those that have two $x_1$'s and two  $y_1$'s.
Lemma~\ref{pinv} tells us that $\tilde g$ is invariant under $T$.
On the other hand, we show that $\tilde h$ is not. It suffices to
show that the sum of all monomials that depend on the
two-dimensional vectors $v^{(1)}$, \dots, $v^{(4)}$ in a
quadrilinear fashion and have type
$\left(\begin{matrix}2\\2\end{matrix}\right)$ is not invariant.
That is clear since the coefficient of the monomial
$x^{(1)}y^{(2)}y^{(3)}y^{(4)}$ will be $1+1+1=1$ after applying
$T$.
\end{proof}

\begin{remark} Replacing the pair $(3,2)$ in the construction of $f$ 
by $(2^t-1,2^t-2)$, where $t$ is an arbitrary natural number, 
we get a multi-linear orthogonal invariant in 
$2((2^t-1)\nu-1)$ vector variables. For  $t>1$, the resulting invariant is 
not a polynomial of the quadratic invariants.
\end{remark}


\section{Two remarks}
\subsection{On the odd-dimensional case}
 As a contrast to  the previous section, we prove the
following theorem. It is a consequence of the first
fundamental theorem for the symplectic group $Sp(2\nu)$, which
holds in its usual form in any characteristic (including 2), as
was proved in \cite[Section 6]{CP}. Using our notation, the first
fundamental theorem  for the symplectic group in characteristic 2
says that the algebra $R_{2\nu\times m}^{Sp(2\nu)}$ is generated
by the $B^{(ij)}$.

\begin{theorem}\label{evenexp}
The invariant $f\in R_{(2\nu +1)\times m}^{O(2\nu +1)}$ is
expressible as a polynomial in the  $Q^{(i)}$ and $B^{(ij)}$ if
and only if the variables $z^{(1)}$, \dots, $z^{(m)}$ occur in $f$
only with even exponents.
\end{theorem}

For example, if $f$ is the square of a (polynomial) invariant,
then $f$ is expressible as a polynomial in the  $Q^{(i)}$ and
$B^{(ij)}$.

\begin{proof}
``Only if" is trivial; we prove ``if".  The proof relies
on the relationship between $O(2\nu +1)$ and $Sp(2\nu)$ that was
described in Subsection~\ref{The-orth}: 
the subalgebra of $R$ generated by the $x$ and $y$ variables is stable with 
respect to the action of $O(2\nu+1)$, and this action can be identified with 
the natural action of $Sp(2\nu)$ on $R_{(2\nu)\times m}$.

Assume hypothesis. 
View $f$ as a polynomial in the variables $z^{(i)}$, and consider a
term $${z^{(1)}}^{2\alpha_1}\cdots
{z^{(m)}}^{2\alpha_m}p\left(x_1^{(1)}, \dots, x_\nu
^{(1)},y_1^{(1)}, \dots, y_\nu ^{(1)}, \dots \dots, x_1^{(m)},
\dots, x_\nu ^{(m)},y_1^{(m)}, \dots, y_\nu ^{(m)}\right)$$ of
highest degree in $f$. Then $p$ must be invariant under
$Sp(2\nu)$, and the first fundamental theorem for the symplectic
group \cite{CP} says that $p$ must be expressible as a polynomial in the
$B^{(ij)}$.

Replace $f$ by the polynomial $$f_1=f-{Q^{(1)}}^{\alpha_1}\cdots
{Q^{(m)}}^{\alpha_m}p.$$ If $f_1$ is expressible in the desired
form, then so is $f$. Of course, $f_1$ is again $O(2\nu
+1)$-invariant, the $z^{(i)}$ occur with even exponents only, and
a highest-degree term of $f$ has disappeared. The  new terms in
$f_1$ are of lower degree. Iterating this procedure, we arrive at
the polynomial 0 after a finite number of steps.
\end{proof}

\subsection{$SO(2\nu)$ versus $O(2\nu)$}
Concerning the even-dimensional case, it is not completely trivial
that $R_{2\nu\times m}^{SO(2\nu)}\neq R_{2\nu\times m}^{O(2\nu)}$
 ($m\geq 2\nu$). An easy proof is possible using a general
theorem of Rosenlicht \cite[Theorem 2]{R} and the fact that
$SO(2\nu)$ is a perfect group (i. e., is generated by commutators
of its elements). We now give an explicit construction of a
$2\nu$-linear polynomial in $R_{2\nu\times 2\nu}$ that is
invariant under $SO(2\nu)$ but not under $O(2\nu)$ \---  just like
the determinant in any characteristic different from 2.

\smallskip

We write $SO(2\nu,\mathbb{C})$ for the special orthogonal group
defined over the complex field by the quadratic form $q=x_1y_1+
\dots +x_\nu y_\nu$. (We continue to write $SO(2\nu)$ for the
group defined over the field $k$ of characteristic 2.) The polar
form $\beta$ of $q$ is given by the same formulas (\ref{betadef}) and (\ref{beta}) 
of Subsection~\ref{The-orth} as over the field $k$.

\begin{lemma}\label{mod2}
If the polynomial $f$ in the coordinates of the indeterminate
$2\nu$-di\-men\-sion\-al vectors $v^{(1)}$, \dots, $v^{(m)}$ has integer
coefficients and is invariant under $SO(2\nu,\mathbb{C})$, then
\---  when viewed as a polynomial over $k$ \---  it is invariant
under $SO(2\nu)$.
\end{lemma}

An analogous statement and proof holds for the groups $O(n,\mathbb{C})$ and
$O(n)$ instead of $SO(2\nu,\mathbb{C})$ and $SO(2\nu)$.

\begin{proof}
For a vector $u\in \mathbb{C}^{2\nu}$ or $u\in k^{2\nu}$,
$q(u)\neq 0$, we write $T_u$ for the reflection in the hyperplane
orthogonal to $u$: $$T_uv\overset{\rm def}=v-\frac{\beta
(v,u)}{q(u)}u.$$

Being invariant under $SO(2\nu,\mathbb{C})$ or $SO(2\nu)$ means
being invariant under the product of any two reflections:
$$f\left(T_uT_wv^{(1)}, \dots, T_uT_wv^{(m)}\right)
=f\left(v^{(1)}, \dots, v^{(m)}\right)$$ for $u,w\in
\mathbb{C}^{2\nu}$  or $u,w\in k^{2\nu}$, $q(u)q(w)\neq 0$.
Coefficients of both sides may be viewed as rational functions
with coefficients in $\mathbb{Z}$ or $\mathbb{Z}/(2)$ of the
vector variables $u$ and $w$, and $SO$-invariance of $f$ boils
down to formal equality of pairs of such rational functions. Since
formal equality over $\mathbb{Z}$ implies that over
$\mathbb{Z}/(2)$, the lemma is proved. 
\end{proof}

We shall use the symbol $*$ to mean any one of the two letters $x$
and $y$.

\begin{proposition}\label{Bszorzat}
Consider the $2\nu$-linear polynomial
 $$\sum
B^{(i_1i_2)}B^{(i_3i_4)}\cdots B^{(i_{2\nu-1}i_{2\nu})}$$ with
integer coefficients, where the $B$'s are defined by
(\ref{alapinv}) over $\mathbb{Z}$, and the sum is extended over
those permutations $i_1$, \dots, $i_{2\nu}$ of the indices 1,
\dots, $2\nu$ that satisfy $i_1<i_2$, $i_3<i_4$, \dots ,
$i_{2\nu-1}<i_{2\nu}$ and $i_1<i_3< \dots <i_{2\nu-1}$.
 The coefficient of the monomial $*_{j_1}^{(1)}\cdots
*_{j_{2\nu}}^{(2\nu)}$ is 1 if $*_{j_1}$, \dots , $*_{j_{2\nu}}$
is a permutation of $x_1$, \dots , $x_\nu$, $y_1$, \dots , $y_\nu$
and is even otherwise.
\end{proposition}

\begin{proof}
The product $$B^{(i_1i_2)}B^{(i_3i_4)}\cdots
B^{(i_{2\nu-1}i_{2\nu})}$$ is the sum of those monomials
$*_{j_1}^{(1)}\cdots *_{j_{2\nu}}^{(2\nu)}$ that satisfy
$j_{i_1}=j_{i_2}$, $j_{i_3}=j_{i_4}$, \dots ,
$j_{i_{2\nu-1}}=j_{i_{2\nu}}$ and have an $x$ and a $y$
corresponding to each of these pairs of indices.
So the sum we are looking at is a linear combination of those
$2\nu$-linear monomials that have the same number \--- say,
$\tau_t$ \---  of $x_t$'s and $y_t$'s among their factors, for
each value of $t$. The coefficient of such a monomial is $\tau_1!
\cdots\tau_\nu!$, since a monomial occurs as many times as its factors
can be grouped into pairs of the form  $\{x_t,y_t\}$. That
coefficient is 1 if  $\tau_1= \dots=\tau_\nu=1$ and even
otherwise.
\end{proof}

Subtract the determinant $D^{(1\cdots (2\nu))}$ from the above sum
(considering both to be  defined over $\mathbb{Z}$). The result is
a polynomial with even coefficients, denote it by $2\Delta$.

\begin{theorem}\label{SOversusO}
The polynomial $\Delta$, viewed as a polynomial over $k$, is invariant
under $SO(2\nu)$ but not under $O(2\nu)$.
\end{theorem}

\begin{proof}
Invariance under $SO(2\nu)$ follows from Lemma~\ref{mod2} as $\Delta$
is invariant under $SO(2\nu,\mathbb{C})$.

Let $*_{j_1}$, \dots , $*_{j_{2\nu}}$ be a permutation of $x_1$,
\dots , $x_\nu$, $y_1$, \dots , $y_\nu$. By
Proposition~\ref{Bszorzat}, the coefficient of the monomial
$*_{j_1}^{(1)}\cdots *_{j_{2\nu}}^{(2\nu)}$ in the polynomial $\Delta$
is 0 if the permutation is even and is  1 if it is odd. It follows
that $\Delta$ is not invariant under the reflection $x_1\leftrightarrow
y_1$ (not even if viewed over $k$), since this transforms the
monomials corresponding to odd permutations into those
corresponding to even ones.
\end{proof}

\section{Separation of orbits}
The  results in this section are analogous to those for
characteristic different from 2. The proofs use Witt's theorem
\cite[Theorem 7.4]{T}, standard facts concerning reductive groups,
and basic algebraic geometry.

Let us introduce the notation  $$A=A_{n\times m}=k\left[Q^{(i)},
B^{(ij)}: 1\leq i\leq m,\; 1\leq i<j\leq m\right].$$ Note that we
have shown in Subsection~\ref{even} that $A\neq R^{O(n)}$ for even
$n$ and large $m$. The same is obvious for odd $n$ and $m\geq n$
as $D^{(1\cdots n)}\in R^{O(n)}\backslash A$.

\subsection{The null-cone}
Recall that the null-cone corresponding to a graded algebra of
polynomials is defined to be the locus of common zeros of its
homogeneous elements of positive degree.

\begin{theorem}\label{null-cone}
The null-cones corresponding to the three algebras
$$R^{SO(n)}_{n\times m} \geq
R^{O(n)}_{n\times m}\geq A_{n\times m}$$ are the same. 
\end{theorem}

\begin{proof}
Suppose that the point $\left(v^{(1)}, \dots,  v^{(m)}\right)$
belongs to the null-cone of $A$; that is, the vectors $v^{(1)}$,
\dots,  $v^{(m)}$ satisfy the equations $Q^{(i)}=0$ and
$B^{(ij)}=0$. The subspace they span is then totally singular (i.
e., has $q\equiv 0$). Let $W$ be a maximal totally singular
subspace containing them. It follows from Witt's theorem that the
dimension of $W$ is $\nu=[n/2]$, and that there exists a maximal
totally singular subspace $W_1$ such that $$k^n=W\oplus
W_1\oplus\ker\beta.$$ For $0\neq t\in k$, let $A_t\in O(n)$ stand
for the special orthogonal transformation that multiplies vectors
in $W$ by $t$, vectors in $W_1$ by $1/t$, and vectors in
$\ker\beta$ by 1. Any $f\in R_{n\times m}^{SO(n)}$ is invariant
under $A_t$, so $$f\left(tv^{(1)}, \dots,
tv^{(m)}\right)=f\left(v^{(1)}, \dots, v^{(m)}\right).$$ This
holds for arbitrary $t\neq 0$, so it must also hold for $t=0$.
This means that the point $\left(v^{(1)}, \dots, v^{(m)}\right)$
is contained in the null-cone of $R^{SO(n)}$.
\end{proof}

\begin{corollary}\label{A-mod}
The algebras $R^{O(n)}$ and $R^{SO(n)}$ are  finitely generated as
$A$-modules.
\end{corollary}

\begin{proof}
Let $G$ stand for $O(n)$ or $SO(n)$. Then $G$ is a reductive
algebraic group, so Nagata's theorem \cite[Theorem 3.4]{N} says
that $R^G$ is finitely generated as an algebra.

Consider a homogeneous element $h\in R^{G}$. By
Theorem~\ref{null-cone} and Nullstellensatz, $h$ has a power in
the ideal of $R$ generated by the $Q^{(i)}$ and the $B^{(ij)}$. It follows by \cite[Lemma
3.4.2]{N} that $h$ has a power in the ideal of $R^G$ generated by
the $Q^{(i)}$ and the $B^{(ij)}$.

 Applying that to each element $h$ of a finite system of homogeneous generators of the
 algebra $R^{G}$ shows that the ideal of $R^G$ generated by the $Q^{(i)}$ and
 the $B^{(ij)}$ contains all elements
 of $R^{G}$ that are homogeneous of high enough degree. So $R^{G}$, as an
 $A$-module, is generated by  elements of degree lower than some number $d$. These form a
 finite-dimensional vector space, so a finite number of them will
 suffice.
\end{proof}

\subsection{Algebro-geometric lemmas}
We recall some well-known facts from algebraic geometry. The word
`variety' below stands for an irreducible affine algebraic variety
over $k$ (the characteristic of $k$ is $2$ in our applications,
but the following general statements are valid if $k$ is an
arbitrary algebraically closed field). Write $K[X]$ for the
algebra of polynomial functions on $X$, and write $K(X)$ for the
field of rational functions on $X$. Let $f:X\to Y$ be a dominant
morphism of varieties. Then the comorphism $f^*$ identifies $K(Y)$
with the subfield $f^*K(Y)$ of $K(X)$. The morphism $f$ is said to
be {\it separable}, if $K(X)\geq f^*K(Y)$ is a separable field
extension. We need the following criterion for separability, see
for example \cite[(17.3) Theorem]{b}: The morphism $f$ is
separable if and only if there is a non-singular point $x$ on $X$
such that $f(x)$ is non-singular in $Y$, and the differential
$d_xf: T_xX\to T_{f(x)}Y$ at $x$ is surjective.

\begin{lemma} \label{lemma-constant-fibres}
Let $f:X\to Y$ be a dominant, separable morphism of varieties.
Suppose that $h$ is a rational function on $X$, such that for some
non-empty Zariski open subset $U$ of $X$, the restriction
$h\vert_U$ is constant along the fibers of $f\vert_U$. Then $h$ is
the pull-back of a rational function on $Y$, that is, $h\in
f^*K(Y)$.
\end{lemma}

\begin{proof} Take a principal affine open subset $V$ in $X$, where
$h\vert_V$ is regular, and $h\vert_V$ is constant along the fibers
of $f\vert_V$. Then $h$ is purely inseparable over $f^*K(Y)$ by
\cite[(18.2) Proposition, p.78]{b}; that is, $h^{p^s}$ is
contained in $f^*K(Y)$ for some natural number $s$. Thus $h$
itself is contained in $f^*K(Y)$, because $f$ is separable by our
assumption.
\end{proof}

More can be said when $Y$ is normal. See for example \cite[(18.3),
p.79]{b}:

\begin{lemma} \label{lemma-normal-image}
Let $f:X\to Y$ be a surjective morphism of varieties, and assume
that $Y$ is normal. Suppose that $h$ is a polynomial function on
$X$, such that $h$ is the pull-back of a rational function on $Y$,
i.e. $h$ is contained in $f^*K(Y)$. Then $h$ is the pull-back of a
polynomial function on $Y$, that is, $h\in f^*K[Y]$.
\end{lemma}

\begin{proof} See for example
\cite[(18.3), p.79]{b}, and note that since we are dealing with
affine varieties, `regular functions' in the sense of \cite{b}
(i.e. everywhere defined rational functions) are the same as
`polynomial functions'.
\end{proof}

\subsection{Rational invariants}
We now look at the field $K^{O(n)}$, which is much easier to deal
with than the algebra $R^{O(n)}$. Note that $K^{O(n)}$ is the
fraction field of $R^{O(n)}$ (this follows easily from the fact
that $SO(n)$ is perfect). 

\begin{theorem}\label{field}
\begin{itemize}
\item[(i)] The field $K_{2\nu\times m}^{O(2\nu)}$ is generated by the
algebraically independent invariants $$Q^{(i)}\qquad\qquad (1\leq
i\leq \min(m,2\nu))$$ and $$B^{(ij)}\qquad\qquad (1\leq i<j\leq
m,\quad i\leq 2\nu).$$

\item[(ii)] For $m<2\nu$ we have 
$K_{2\nu\times m}^{SO(2\nu)}
=K_{2\nu\times m}^{O(2\nu)}$. 
For $m\geq 2\nu$, the field $K_{2\nu\times m}^{SO(2\nu)}$ is a quadratic 
extension of $K_{2\nu\times m}^{O(2\nu)}$, generated for example by the 
invariant $\Delta$ constructed in Theorem~\ref{SOversusO}. 

\item[(iii)]  The field $K_{(2\nu+1)\times m}^{SO(2\nu+1)}$ is generated
by the algebraically independent invariants $$Q^{(i)}\qquad\qquad
(1\leq i\leq \min(m,2\nu)),$$ $$B^{(ij)}\qquad \qquad (1\leq
i<j\leq m,\quad i\leq 2\nu),$$ and $$D^{(1, \dots,
2\nu,l)}\qquad\qquad (2\nu +1\leq l\leq m).$$
\end{itemize}
\end{theorem}

The description in 
(ii) of $K_{2\nu\times m}^{SO(2\nu)}$  for $m\geq  2\nu$  will be made 
complete in  Theorem~\ref{gamma} where we determine the quadratic 
polynomial over $K_{2\nu\times m}^{O(2\nu)}$ that $\Delta$  satisfies.

Note that the theorem is valid in any characteristic. In any
characteristic different from 2, 
the third statement remains valid if $SO(2\nu+1)$ is replaced by
$O(2\nu+1)$ and $D^{(1, \dots, 2\nu,l)}$ is  replaced by
$B^{(2\nu+1\mid l)}$. The proof given below, appropriately
modified, goes through.

The proof is via   the following propositions.

\begin{proposition}\label{betaszurj}
 Let $m$ be any positive integer, and let $\left(\beta
^{(ij)}\right)$ be any alternating  $m\times m$ matrix of rank $r\leq n$.
Then there exist vectors $u^{(1)}, \dots , u^{(m)}\in k^{n}$
with $$\beta\left(u^{(i)},u^{(j)}\right)=\beta^{(ij)} \qquad
(i,j=1, \dots, m).$$
\end{proposition}

\begin{proof}
It is well known that $\left(\beta ^{(ij)}\right)$ is cogredient
to $J\oplus 0= \left(
\begin{matrix} 0 & I \\ I & 0 \end{matrix} \right) \oplus 0$
with $J$ of size $r\times r$ (so  $r$ is always even). The
proposition obviously holds for the latter matrix, and the general
case follows by base change.
\end{proof}

\begin{proposition}\label{bqszurj}
 Let  $m\leq n$. Let $\left(\beta ^{(ij)}\right)$ be any
$m\times m$ alternating matrix, and let $q^{(1)}, \dots ,
q^{(m)}\in k$. Then there exist vectors $v^{(1)}, \dots ,
v^{(m)}\in k^n$ with
$$\beta\left(v^{(i)},v^{(j)}\right)=\beta^{(ij)} \qquad (i,j=1,
\dots, m) $$ and $$q\left(v^{(i)}\right)=q^{(i)} \qquad (i=1,
\dots , m).$$
\end{proposition}

\begin{proof}
As always, we set $\nu=[n/2]$. Choose vectors $u^{(1)}, \dots ,
u^{(m)}\in k^{2\nu}$ as in the previous proposition.

Consider $n=2\nu +1$ first. 
Note that the standard quadratic form $q$ is onto $k$ on any line
parallel to $\ker\beta$ (the $z$-axis). 
Therefore, there exist vectors  $v^{(i)}\in k^{2\nu +1}$ that are
mapped  to the $u^{(i)}$  by the projection $$k^{2\nu +1}\to
k^{2\nu +1}/\ker \beta =k^{2\nu} $$ and have
$q\left(v^{(i)}\right)=q^{(i)}$. 

Now let $n=2\nu$. First suppose that $m=n$ and $u^{(1)}, \dots ,
u^{(m)}$ is a basis of $k^n$. Define a new quadratic form $q^*$ by
the formula $$ q^* \left(\sum_{i=1}^m \lambda_i
u^{(i)}\right)=\sum_{i=1}^m\lambda_i^2q^{(i)} + \sum_{1\leq i <j
\leq m} \lambda_i \lambda_j \beta ^{(ij)}.$$  
Let $\beta^*$ stand
for the polar form of $q^*$. Then $$ \beta ^*
\left(u^{(i)},u^{(j)}\right)=q^*\left(u^{(i)}+u^{(j)}\right)-
q^*\left(u^{(i)}\right)-q^*\left(u^{(j)}\right)=\beta^{(ij)}=\beta
\left(u^{(i)},u^{(j)}\right),$$ therefore, $\beta^*\equiv \beta$.
It follows that  $q^*$ is non-degenerate. Since $k$ is
algebraically closed, all non-degenerate quadratic forms are
equivalent. So there is a linear isomorphism $A:k^n \to k^n$ such
that $q(Au)=q^*(u)$ for all $u\in k^n$. It of course follows that
$$\beta \left(Au', Au''\right)=\beta ^* \left(u', u''\right)=\beta
\left(u', u''\right)$$ for all  $u', u''\in k^n$ (that is, $A\in
Sp(n)$). Define $v^{(i)}=Au^{(i)}$ \;$ (i=1, \dots , m)$. Then $$
\beta \left(v^{(i)}, v^{(j)}\right)= \beta \left(u^{(i)},
u^{(j)}\right)=\beta^{(ij)}$$ and $$
q\left(v^{(i)}\right)=q^*\left(u^{(i)}\right)=q^{(i)},$$ i. e.,
$v^{(1)}$, \dots ,  $v^{(m)}$ have the desired properties.

Suppose finally that $n=2\nu$ but $u^{(1)}, \dots , u^{(m)}$ do
not span $k^n$. Choose some vector $0\neq u^{(0)}\in \left\langle
u^{(1)}, \dots , u^{(m)}\right\rangle ^\bot$. Choose a linear
function $f:k^n \to k$ with $f\left(u^{(0)}\right)\neq 0$.
 Define the new quadratic form $q^*$ by the formula $$ q^* =q + \lambda f^2,$$ with some $\lambda
\in k$ that gives $q^*\left(u^{(0)}\right)\neq 0$. The quadratic
form $f^2$ has 0 as its polar form, so $q^*$ has $\beta$. It
follows that $q^*$ is non-degenerate. We therefore have a linear
isomorphism $A:k^n \to k^n$ such that $q(Au)=q^*(u)$ for all $u\in
k^n$. Of course $A\in Sp(n)$. The vectors $Au^{(i)}$ have
$$\beta\left(Au^{(i)}, Au^{(j)}\right)=\beta \left(u^{(i)},
u^{(j)}\right)=\beta ^{(ij)}.$$ Note also that
$Au^{(0)}\in\left\langle Au^{(1)}, \dots ,  Au^{(m)}\right\rangle
^\bot$ and $q\left(Au^{(0)}\right)\neq 0$. The latter ensures that
$q$ is onto $k$ on any line parallel to $kAu^{(0)}$. So there are
vectors $v^{(i)}\in Au^{(i)}+kAu^{(0)}$ with
$q\left(v^{(i)}\right)=q^{(i)}$. They have all desired properties.
\end{proof}

We shall need the following consequence of Witt's theorem.

\begin{proposition}\label{Wittalt}
\begin{itemize} 
\item[(i)] For $n=2\nu$  and arbitrary $m$, there exists a non-empty open
set $U\subset k^{n\times m}$ with the following property: if
$$\left({v^{(1)}} ', \dots, {v^{(m)}}'\right)\in U \textrm{ and }
\left({v^{(1)}}'', \dots, {v^{(m)}}''\right)\in U$$ satisfy
$$\begin{aligned} {Q^{(i)}}'&={Q^{(i)}}'' &(1\leq i\leq 2\nu),\\
{B^{(ij)}}' & ={B^{(ij)}}''\qquad &(1\leq i<j\leq m,\quad i\leq
2\nu),\end{aligned}$$  then there is an orthogonal transformation
$A$ such that $A{v^{(i)}}'={v^{(i)}}''$ for every $1\leq i\leq m$.

\item[(ii)] When $m<n=2\nu$, the assertion (i) holds with $A$ taken from 
the special orthogonal group $SO(2\nu)$. 

\item[(iii)] For $n=2\nu+1$  and arbitrary $m$, there exists a non-empty
open set $U\subset k^{n\times m}$ with the following property: if
$$\left({v^{(1)}} ', \dots, {v^{(m)}}'\right)\in U {\textrm{ and
}} \left({v^{(1)}}'', \dots, {v^{(m)}}''\right)\in U$$ satisfy
$$\begin{aligned} {Q^{(i)}}' & ={Q^{(i)}}''\qquad \qquad \qquad
&(1\leq i\leq 2\nu) ,\\ {B^{(ij)}}' & ={B^{(ij)}}''\qquad &(1\leq
i<j\leq m,\quad i\leq 2\nu)  ,\\ {D^{(1,\cdots,2\nu,l)}}' &
={D^{(1,\cdots,2\nu,l)}}''\quad\qquad\qquad &( 2\nu+1\leq l\leq m)
 ,\end{aligned} $$  then there is an orthogonal transformation $A$ such that $A{v^{(i)}}'={v^{(i)}}''$ for every
$1\leq i\leq m$.
\end{itemize} 
\end{proposition}

\begin{proof}
For (i) and (iii), 
an $m$-tuple of vectors shall be contained in $U$ exactly if the
images of the  first $\min (m,2\nu)$ vectors are linearly
independent in $k^n/\ker\beta$. Those first $\min (m,2\nu)$
vectors will always span a  subspace $W$ with $W\cap\ker\beta=0$.
If $$\left({v^{(1)}} ', \dots, {v^{(m)}}'\right)\in U {\textrm{
and }} \left({v^{(1)}}'', \dots, {v^{(m)}}''\right)\in U$$ satisfy
the conditions stated in the proposition, then Witt's theorem
provides
 $A\in O(n)$ with
$$A{v^{(i)}}'={v^{(i)}}''\qquad (i=1, \dots, \min(m,2\nu)).$$ If
$m>2\nu$, we need to show that this equality also holds  for
$2\nu<i\leq m$.

\smallskip

(i) As $\beta $ is non-degenerate and ${v^{(1)}}''$, \dots,
${v^{(2\nu)}}''$ is a basis of $k^{2\nu}$, it suffices to show
that
\begin{equation}\label{odamegy}
\beta\left({v^{(i)}}'',A{v^{(j)}}'\right)=
\beta\left({v^{(i)}}'',{v^{(j)}}''\right)
\end{equation} for $1\leq j\leq m$ and $1\leq i\leq 2\nu$.
 This is equivalent to
 $$\beta\left(A{v^{(i)}}',A{v^{(j)}}'\right)=\beta\left({v^{(i)}}',{v^{(j)}}'\right),$$
 which follows from the orthogonality of $A$.

\smallskip

(iii) Equality (\ref{odamegy}) is proved as above, and shows that
$A{v^{(j)}}'$ and ${v^{(j)}}''$ can differ only in their $z$
coordinates.  Equality of the $z$ coordinates will follow from
\begin{equation}\label{det=}\det\left[{v^{(1)}}'', \dots, {v^{(2\nu)}}'',A{v^{(j)}}'\right]=
\det\left[{v^{(1)}}'', \dots,
{v^{(2\nu)}}'',{v^{(j)}}''\right],\end{equation} since expanding
both determinants by the last column gives the same terms except
for the term containing the $z$ coordinate of the last vector with
the same non-vanishing $2\nu\times 2\nu$ minor as its coefficient
on both sides.

Equality (\ref{det=}) is equivalent to $$\det\left[A{v^{(1)}}',
\dots, A{v^{(2\nu)}}',A{v^{(j)}}'\right] =\det\left[{v^{(1)}}',
\dots, {v^{(2\nu)}}',{v^{(j)}}'\right],$$ which follows from
orthogonality of $A$.

\smallskip 

(ii) We impose an additional condition on points of $U$: 
the orthogonal subspace to the subspace spanned by the 
components of an  $m$-tuple ($m<2\nu$) in $U$ should 
contain a non-singular vector. It is easy to see that $U$ still 
contains a non-empty Zariski 
open subset in $k^{2\nu\times m}$. 
Indeed, when $m=2\nu-1$, the orthogonal subspace to the subspace spanned by 
the linearly independent components of an $m$-tuple $v\in k^{2\nu\times m}$ 
is spanned by a vector whose coordinates are $m\times m$ minors of $v$, 
therefore the condition that this vector is non-singular is expressed as 
the non-vanishing of a polynomial function on $k^{2\nu\times m}$. 
($U$ is clearly non-empty; for example, a basis of the subspace orthogonal 
to some non-singular vector is contained in $U$.) 
To handle the case $m<2\nu-1$ as well, note that the image of 
a non-empty Zariski open subset of $k^{2\nu\times(2\nu-1)}$ 
under  the projection map onto $k^{2\nu\times m}$ contains a 
non-empty open subset of $k^{2\nu\times m}$. 

Now take from $U$ the $m$-tuples $v'$, $v''$ satisfying the 
conditions stated in the proposition. 
By Witt's theorem we have $A\in O(2\nu)$ with $Av'=v''$. 
There is a non-singular 
vector $u$ orthogonal to the subspace spanned by the 
components of $v'$. The reflection $T_u$ fixes 
$v'$. So both $A$ and $AT_u$ map $v'$ to $v''$, 
and one of them is contained in $SO(2\nu)$.

\end{proof}

\begin{proofof}{Theorem~\ref{field}} 
(i) and (iii): Write 
$f$ for the regular map defined on $k^{n\times m}$ that has the
invariants in the theorem as its coordinates. For $m\leq 2\nu$,
Proposition~\ref{bqszurj} shows that $f$ is surjective. If $m\geq
2\nu$, $f$ is still dominant, for if we prescribe values
$q^{(i)}$, $\beta^{(ij)}$ and (in the odd-dimensional case)
$d^{(1,\dots,2\nu,l)}$ with
$\det\left(\beta^{(ij)}\right)_{i,j=1}^{2\nu}\neq 0$, then the
vectors  $v^{(1)}$, \dots, $v^{(2\nu)}$ provided by
Proposition~\ref{bqszurj} will give a basis in
$k^n/\ker\beta=k^{2\nu}$, and this ensures the existence of
$v^{(2\nu+1)}$, \dots , $v^{(m)}$ such that the coordinates of $f$
take the prescribed values on the  $m$-tuple $v^{(1)}$, \dots,
$v^{(m)}$. This proves algebraic independence of the invariants in
the theorem.

We now show that $f$ is separable. Consider the  point
$\left(e^{(1)}, \dots, e^{(m)}\right)$ in $k^{n\times m}$ given by the
first $\min(m,2\nu)$ vectors of the standard basis of $k^n$ and
$m-\min(m,2\nu)$ zero vectors. We claim that the differential of
$f$ at this point is onto. The partial derivatives are as follows.
\begin{equation}\label{Qparc}
\frac{\partial Q^{(i)}}{\partial x_t^{(i)}}=y_t^{(i)},\qquad
\frac{\partial Q^{(i)}}{\partial y_t^{(i)}}=x_t^{(i)},
\end{equation}
all other partials of $Q^{(i)}$ being zero. So the $n\times m$
matrix formed by the partials of $Q^{(i)}$ has $e^{(i)}$ with $x$ and
$y$ coordinates  interchanged as its $i$th column, all other
columns being zero. Also,
\begin{equation}\label{Bparc}
\frac{\partial B^{(ij)}}{\partial x_t^{(i)}}=y_t^{(j)},\qquad
\frac{\partial B^{(ij)}}{\partial y_t^{(i)}}=x_t^{(j)},\qquad
\frac{\partial B^{(ij)}}{\partial x_t^{(j)}}=y_t^{(i)},\qquad
\frac{\partial B^{(ij)}}{\partial y_t^{(j)}}=x_t^{(i)},
\end{equation}
all other partials of $B^{(ij)}$ being zero.  So the $n\times m$
matrix formed by the partials of $B^{(ij)}$ has  $e^{(i)}$ with $x$
and $y$ coordinates  interchanged as its $j$th column and has
$e^{(j)}$ with $x$ and $y$ coordinates  interchanged as its $i$th
column, all other columns being zero. We easily see that all these
$n\times m$ matrices are linearly independent. Our claim follows
in the even-dimensional case; in the odd-dimensional case we
observe that these $(2\nu+1)\times m$ matrices  have nothing but
zeros in their last lines, so it suffices to prove that the last
lines of the $(2\nu+1)\times m$ matrices formed by the partials of the
$D^{(1,\dots, 2\nu,l)}$ are linearly independent. This is obvious,
since $$\frac{\partial D^{(1, \dots, 2\nu,l)}}{\partial
z^{(l')}}=\delta_{(l')}^{(l)}$$ for  $2\nu +1\leq l,l'\leq m$.

Now let $h\in K_{n\times m}^{O(n)}$ (considered as a function on
$k^{n\times m}$). 
Then $h$ is constant along the orbits of $O(n)$,
so Proposition~\ref{Wittalt} shows that $h$ is constant along the
fibers of $f$ (at least on some non-empty open set).  By
Lemma~\ref{lemma-constant-fibres}, $h$ is the  pull-back of a
rational function. 

(ii) When $m<2\nu$, the same argument as above works: 
$h\in K_{2\nu\times m}^{SO(2\nu)}$ is constant along the fibers of $f$ 
defined above by Proposition~\ref{Wittalt} (ii), 
so by Lemma~\ref{lemma-constant-fibres}, $h$ is 
a rational function in the $Q^{(i)}$, $B^{(ij)}$. 

For the case $m\geq 2\nu$, note that $K^{O(2\nu)}$ is the fixed point set 
of the two-element group $O(2\nu)/SO(2\nu)$ acting on $K^{SO(2\nu)}$, 
hence the degree of the field extension $K^{SO(2\nu)}\mid K^{O(2\nu)}$ 
is $1$ or $2$. By Theorem~\ref{SOversusO}, it must be a quadratic extension 
generated by $\Delta$. 

\end{proofof}

\subsection{The case $m\leq n$}
The results of this section show that the conjectured exotic orthogonal
invariants can appear only if the number of vector variables is sufficiently
large, namely, if  $m>n$.

\begin{theorem}\label{matmost2nu}
Let $n=2\nu$ or $n=2\nu +1$, and let $m\leq 2\nu$. Then the
algebra
$R_{n\times m}^{O(n)}$ 
is generated by the $\binom{m+1}{2}$ algebraically independent
invariants  $Q^{(i)}$ and $B^{(ij)}$. 
When $m<2\nu=n$, we have 
$R_{2\nu\times m}^{SO(2\nu)}=R_{2\nu\times m}^{O(2\nu)}$. 
\end{theorem}

\begin{proof}
Let $f:k^{n\times m}\to k^{\binom{m+1}{2}}$ stand for the regular
map that has the $Q^{(i)}$ and $B^{(ij)}$ as its coordinates.
Choose any $h\in R_{n\times m}^{O(n)}$ 
(or $h\in R_{2\nu\times m}^{SO(2\nu)}$ when $m<2\nu$). 
Theorem~\ref{field} says
that $h$ is the pull-back of a rational function. But $h$ is a
polynomial, and Proposition~\ref{bqszurj} says that $f$ is
surjective. By Lemma~\ref{lemma-normal-image}, $h$ is the
pull-back of a polynomial function.
\end{proof}

\newcommand{\Z}{\mathbb{Z}}
\newcommand{\Q}{\mathbb{Q}}

\begin{theorem}\label{R}
 Let $m=n=2\nu+1$. Let $D$ stand for $D^{(1\cdots n)}$. Then the
 algebra
$R_{n\times m}^{O(n)}$ 
is generated by the 
$\binom{n+1}{2}+1$ invariants $Q^{(i)}$, $B^{(ij)}$ and $D$, the ideal
of algebraic relations between whom  is generated by the single
element $G$ defined as $$G=D^2-\frac{1}{2}\left|
\begin{matrix}2Q^{(1)} & B^{(12)} &\cdots & B^{(1n)}\\B^{(21)} &
2Q^{(2)} & \cdots & B^{(2n)}\\\vdots & \vdots & \ddots & \vdots\\
B^{(n1)} & B^{(n2)} & \cdots & 2Q^{(n)}\end{matrix}\right|.$$
(See Proposition~\ref{Gwelldef} for the meaning of 1/2 here.)
\end{theorem}

We break the proof up into several propositions.

\begin{proposition}\label{Gwelldef}
The determinant in the definition of $G$, when interpreted 
as a polynomial over
$\Z$ in the variables $Q^{(i)}$ and $B^{(ij)}$, has even coefficients. 
So $G$ is defined as a polynomial
over $\Z$ and {\rm a fortiori} over $k$.
\end{proposition}

\begin{proof}
Each expansion term in the determinant either has a factor from
the diagonal and therefore has an even coefficient, or is a
product of off-diagonal entries and can be paired with the
transposed term (note that $B^{(ij)}=B^{(ji)}$).
\end{proof}

\begin{proposition}\label{G=0}
The polynomials $Q^{(i)}$, $B^{(ij)}$ and $D$ satisfy the relation
$$(-1)^\nu D^2-\frac{1}{2}\left|
\begin{matrix}2Q^{(1)} & B^{(12)} &\cdots & B^{(1n)}\\B^{(21)} &
2Q^{(2)} & \cdots & B^{(2n)}\\\vdots & \vdots & \ddots & \vdots\\
B^{(n1)} & B^{(n2)} & \cdots & 2Q^{(n)}\end{matrix}\right|=0$$
over $\Z$ and {\rm a fortiori} over $k$.
\end{proposition}

\begin{proof}
Working over $\Q$, the matrix of the polar form $\beta$ of the
quadratic form $$q=x_1y_1+\dots+x_\nu y_\nu +z^2$$ is
$$M
=\left(\begin{matrix}0&1\\1&0\end{matrix}\right)\oplus\dots\oplus\left(\begin{matrix}0&1\\1&0\end{matrix}\right)\oplus(2).$$
For arbitrary $V\in \Q^{n\times n}$ with $i$th column $v^{(i)}$, we
have $$V^TMV=\left(\beta\left(v^{(i)},
v^{(j)}\right)\right)_{i,j=1}^n.$$ Taking determinants gives
$$(-1)^\nu\cdot 2 \cdot (\det V)^2=\det\left(\beta\left(v^{(i)},
v^{(j)}\right)\right)_{i,j=1}^n.$$ The proposition follows, since
$\beta \left( v^{(i)}, v^{(i)}\right)=2q\left(v^{(i)}\right)$.
\end{proof}

The following proposition deals with the hypersurface $\{G=0\}$ 
in the affine space $k^{\binom{n+1}{2}+1}$, with 
coordinates denoted by $Q^{(i)}$, $B^{(jl)}$, $D$ 
($1\leq i\leq n$, $1\leq j<l\leq n$). 

\begin{proposition}\label{prop:normal}
The  hypersurface $\{G=0\}$ in
$k^{\binom{n+1}{2}+1}$ is normal.
\end{proposition}

\begin{proof}
A hypersurface $H$ (the zero locus of a single polynomial in an
affine space) is normal if and only if the set of singular points
has codimension $\geq 2$ in $H$; this follows for example from
Seidenberg's criterion for normality \cite[Theorem 3]{S}, together
with Macaulay's unmixedness theorem (cf. \cite[Theorem 17.6]{m}).

Write $G$ as $$G=D^2-\left(Q^{(1)}
F^{(1)}+\dots+Q^{(n)}F^{(n)}+F^{(0)}\right),$$ where $F^{(i)}$ is
the $i$th principal $(n-1)\times (n-1)$ minor of the matrix
$$\left(\begin{matrix}2Q^{(1)} & B^{(12)} &\cdots &
B^{(1n)}\\B^{(21)} & 2Q^{(2)} & \cdots & B^{(2n)}\\\vdots & \vdots
& \ddots & \vdots\\ B^{(n1)} & B^{(n2)} & \cdots &
2Q^{(n)}\end{matrix}\right),$$ and  $F^{(0)}$ is the sum of those
terms in the determinant of this matrix that have no factor from the
diagonal, and have at least $\nu+1$
factors from above the diagonal.

In particular, $$\frac{\partial G}{\partial Q^{(i)}}
=F^{(i)}\qquad(i=1,\dots,n).$$ We claim that the locus of common
zeros of $G$, $F^{(1)}$ and $F^{(n)}$ is of codimension 3 in
$k^{\binom{n+1}{2}+1}$. Equivalently, the locus of common zeros of
$F^{(1)}$ and $F^{(n)}$ is of codimension 2 in the hyperplane
$\{D=0\}$. 
(To see the equivalence note that projection from the direction of the $D$ 
coordinate axis onto the coordinate hyperplane $\{D=0\}$ maps 
the hypersurface $\{G=0\}$ bijectively onto the hyperplane $\{D=0\}$.) 
The polynomials $F^{(1)}$ and $F^{(n)}$ depend only on the variables 
$B^{(ij)}$, and their vanishing on a common hypersurface in the affine 
space $\{D=0\}$ would mean having the defining polynomial of that 
hypersurface as a common factor. 
Therefore it suffices to show that $F^{(1)}$ and $F^{(n)}$ have
no common factors as polynomials in the $B^{(ij)}$. To this end, we
 impose the
order
$$\begin{aligned}B^{(12)}&>&B^{(23)}&>&\cdots\cdots&>&B^{(n-2|n-1)}&>&B^{(n-1|n)}&>\\
                     &>&B^{(13)}&>&\cdots\cdots&>&B^{(n-3|n-1)}&>&B^{(n-2|n)}&>\\
                     &  &       & &\cdots\cdots&& \cdots\cdots       && \cdots\cdots& \\
                     &  &       &  &     & >&B^{(1|n-1)}&>&B^{(2|n)}&>\\
                     &  &       &  &      &&
                     &>&B^{(1|n)}\end{aligned}$$  on the variables and the corresponding lexicographic order on the monomials. Then the leading
                     monomial of $F^{(1)}$ is
                      $\left({B^{(23)}}{B^{(45)}}\cdots{B^{(n-1|n)}}\right)^2$,
                     and the leading monomial of $F^{(n)}$ is
                     $\left({B^{(12)}}{B^{(34)}}\cdots{B^{(n-2|n-1)}}\right)^2$.
                     The leading monomials have no common factors,
                     hence $F^{(1)}$ and $F^{(n)}$ have no common
                     factors. So the locus of common zeros of $G$, $F^{(1)}$ and $F^{(n)}$ is of codimension 3.
                     The singular locus of $\{G=0\}$ is contained
                     in that locus, so it has codimension $\geq 2$ in
$\{G=0\}$, which is therefore normal.
\end{proof}

\begin{proofof}{Theorem~\ref{R}}
Consider the map $$f:k^{n\times n}\to k^{\binom{n+1}{2}+1}$$ that
has the $Q^{(i)}$, the $B^{(ij)}$, and $D$ as its coordinates. It
follows from Propositions~\ref{G=0} and \ref{bqszurj} that the
image of $f$ is the hypersurface $\{G=0\}$ (we need that the
characteristic is 2, so the values of the $Q^{(i)}$ and the $B^{(ij)}$
determine the value of $D$ on $\{G=0\}$).

Choose any $h\in R_{n\times n}^{O(n)}$. Theorem~\ref{field} says
that $h$ is the pull-back of a rational function on $\{G=0\}$. But
$h$ is a polynomial, so, by Lemma~\ref{lemma-normal-image} 
and Proposition~\ref{prop:normal}, $h$ is
the pull-back of a polynomial.
\end{proofof}

\bigskip

We now turn to the description of the algebra  of special orthogonal
invariants in the case $m=n=2\nu$.  We shall write $\sum BB\cdots B$ for the
$2\nu$-linear $O(2\nu,\mathbb C)$-invariant 
$$\sum
B^{(i_1i_2)}B^{(i_3i_4)}\cdots B^{(i_{2\nu-1}i_{2\nu})}$$ defined over $\Z$
that was proved in Proposition~\ref{Bszorzat} to agree with
$D=D^{(1\cdots(2\nu))}$  modulo 2.

\begin{theorem}\label{gamma}
 Let $m=n=2\nu$. Let $\Delta$ stand for the $SO(2\nu)$-invariant constructed in
 Theorem~\ref{SOversusO}.  Then the
 algebra
$R_{n\times m}^{SO(n)}$ 
is generated by the 
$\binom{n+1}{2}+1$ invariants $Q^{(i)}$, $B^{(ij)}$ and $\Delta$, the ideal
of algebraic relations between whom  is generated by the single element 
$\Gamma$ defined as $$\Gamma=\Delta^2-\Delta\sum
BB\cdots B+\frac 14\left(\left(\sum
BB\cdots B\right)^2-(-1)^\nu
\left|
\begin{matrix}2Q^{(1)} & B^{(12)} &\cdots & B^{(1n)}\\B^{(21)} &
2Q^{(2)} & \cdots & B^{(2n)}\\\vdots & \vdots & \ddots & \vdots\\
B^{(n1)} & B^{(n2)} & \cdots & 2Q^{(n)}\end{matrix}\right|\right).$$ 
(See the proof for the meaning of 1/4 here.)
\end{theorem}

\begin{proof} The proof is rather similar to that of Theorem~\ref{R}. 
Write $L$ for the expression 
$$\left(\sum
BB\cdots B\right)^2-(-1)^\nu
\left|
\begin{matrix}2Q^{(1)} & B^{(12)} &\cdots & B^{(1n)}\\B^{(21)} &
2Q^{(2)} & \cdots & B^{(2n)}\\\vdots & \vdots & \ddots & \vdots\\
B^{(n1)} & B^{(n2)} & \cdots & 2Q^{(n)}\end{matrix}\right|,$$ 
so $L$ is a polynomial of $Q^{(i)}$, $B^{(ij)}$ 
with integral coefficients. 

First interpret $Q^{(i)}$, $B^{(ij)}$, $\Delta$ as polynomials over $\Z$ 
in the $x,y$ variables. 
Recall that $\Delta$ was defined over $\Z$ by
$$\Delta=\frac12\left(\sum BB\cdots B -D\right),$$ where
$D=D^{(1\cdots(2\nu))}$. 
Set $$\bar\Delta=\frac12\left(\sum
BB\cdots B +D\right).$$ 
It is a polynomial with integral coefficients in the $x,y$ variables by 
Proposition~\ref{Bszorzat}, hence so is  
$$\Delta\bar\Delta=\frac 14\left(\left(\sum BB\cdots
B\right)^2-D^2\right)=\frac 14 L$$ 
(the second equality is proved in the same manner as 
Proposition~\ref{G=0}).  
Note that 
$\Delta+\bar\Delta=\sum BB\cdots B$. It follows that 
$Q^{(i)}$, $B^{(ij)}$, $\Delta$ (considered as polynomials over $\Z$ 
in the $x,y$ variables) satisfy the relation 
\begin{equation}\label{eq:integralgamma} 
\Delta^2-\Delta\sum BB\cdots B + L/4=0.
\end{equation} 
We claim that the coefficients of $L$ are divisible by four, 
so $L/4$ is a polynomial in the variables $Q^{(i)}$, $B^{(ij)}$ with 
integer coefficients. Indeed, multiply the relation 
$\Delta\bar\Delta=L/4$ by $4$ and consider it modulo $2$:   
the left hand side becomes zero, 
so we obtain on the right hand side an algebraic relation over $k$ holding 
between $Q^{(i)}$, $B^{(ij)}$ (defined over $k$). But 
$Q^{(i)}$, $B^{(ij)}$ are algebraically independent in 
$R_{n\times m}^{SO(n)}$ by Theorem~\ref{field}, 
so this relation must be trivial. This means that all coefficients 
of $L$ (as a polynomial in the $Q^{(i)}$, $B^{(ij)}$) are even. 
Taking now the relation $2\Delta\bar\Delta=L/2$ modulo $2$ 
and repeating the same argument we obtain our claim.   
So \eqref{eq:integralgamma} is an algebraic relation with 
integral coefficients holding between $Q^{(i)}$, $B^{(ij)}$, $\Delta$ 
(considered as polynomials over $\Z$ in the $x,y$ variables). 

It follows immediately that \eqref{eq:integralgamma}  
makes sense and holds as a relation over $k$; that is, 
the relation $\Gamma=0$ makes sense and holds in $R_{n\times m}^{SO(n)}$. 

Consider now the map $$f:k^{n\times n}\to k^{\binom{n+1}{2}+1}$$ that
has the $Q^{(i)}$, the $B^{(ij)}$, and $\Delta$ as its coordinates. It
follows from the relation $\Gamma=0$ and Proposition~\ref{bqszurj} that the
image of $f$ is the hypersurface $\{\Gamma=0\}$ in 
$k^{\binom{n+1}{2}+1}$. 
(For surjectivity, we also need that the coset $O(2\nu)\backslash
SO(2\nu)$ interchanges $\Delta$ and $\bar\Delta$, so a point
$(Q,B,\Delta)$ is in the image of $f$ 
if and only if 
$(Q,B,\bar\Delta)$ is in the image of $f$.) 
Choose any $h\in R_{n\times n}^{SO(n)}$. Theorem~\ref{field} says
that $h$ is the pull-back of a rational function on the hypersurface 
$\{\Gamma=0\}$. But $h$ is a polynomial, so, by
Lemma~\ref{lemma-normal-image} and Proposition~\ref{prop:gammanormal} below, 
$h$ is the pull-back of a polynomial.
\end{proof} 

\begin{proposition}\label{prop:gammanormal} 
Consider the affine space $k^{\binom{n+1}2 +1}$, with
coordinates denoted by 
$Q^{(i)}$, $B^{(jl)}$, $\Delta$ ($1\leq i\leq n$, $1\leq j<l\leq n$). 
Then the hypersurface $\{\Gamma=0\}$ in $k^{\binom{n+1}{2}+1}$ is
normal. 
\end{proposition}

\begin{proof}
Just as in Proposition~\ref{prop:normal}, it suffices to prove that the
singular locus has codimension $\geq 2$ in the hypersurface.

 Calculate $$\frac{\partial\Gamma}{\partial\Delta}=\sum BB\cdots
 B$$ and  $$\frac{\partial\Gamma}{\partial Q^{(n)}}=\frac12\left|\begin{matrix}2Q^{(1)} & B^{(12)} &\cdots &
B^{(1|n-1)}\\B^{(21)} & 2Q^{(2)} & \cdots & B^{(2|n-1)}\\\vdots &
\vdots & \ddots & \vdots\\ B^{(n-1|1)} & B^{(n-1|2)} & \cdots &
2Q^{(n-1)}\end{matrix}\right|=Q^{(1)}
F^{(1)}+\dots+Q^{(n-1)}F^{(n-1)}
+F^{(0)},$$ where $F^{(i)}$ is the $i$th
principal $(n-2)\times (n-2)$ minor of the last determinant 
for $i=1,\ldots,n-1$, and $F^{(0)}$ also depends only on the $B^{(ij)}$. 

 We claim that the locus of common
zeros of $\Gamma$, $\partial\Gamma/\partial\Delta$ and
$\partial\Gamma/\partial Q^{(n)}$ is of codimension 3 in
$k^{\binom{n+1}{2}+1}$. Equivalently, the locus of common zeros of
$\partial\Gamma/\partial\Delta$ and $\partial\Gamma/\partial
Q^{(n)}$  is of codimension 2 in the hyperplane $\{\Delta=0\}$. It
suffices to show that $\partial\Gamma/\partial\Delta$ and
$\partial\Gamma/\partial Q^{(n)}$  have no common factors as
polynomials in the $Q^{(i)}$ and the $B^{(ij)}$. As
$\partial\Gamma/\partial\Delta$ depends only on the $B^{(ij)}$, so
will any common factor, but then, in order to divide
$\partial\Gamma/\partial Q^{(n)}$, it must divide each $F^{(i)}$.
But we have shown in the proof of Proposition~\ref{prop:normal} that
$F^{(1)}$ and $F^{(n-1)}$   (there denoted by $F^{(1)} $ and
$F^{(n)}$ since $n$ there was odd and the $F$'s were $(n-1)\times
(n-1)$ minors of an $n\times n$ matrix) have no common factors. So
the locus of common zeros of $\Gamma$,
$\partial\Gamma/\partial\Delta$ and $\partial\Gamma/\partial
Q^{(n)}$  is of codimension 3.
                     The singular locus of $\{\Gamma=0\}$ is contained
                     in that locus, so it has codimension $\geq 2$ in $\{\Gamma=0\}$, which therefore is normal.
\end{proof}

\medskip
\noindent{\bf Acknowledgment.} 
We thank the referee for asking about rational $SO(2\nu)$-invariants. 
This led to the statement (ii) in Theorem~\ref{field}, 
and to a description of the polynomial 
$SO(2\nu)$-invariants for $m\leq 2\nu$.


\begin{thebibliography}{99999}
\bibitem{b} A. Borel,
Linear Algebraic Groups, W. A. Benjamin Inc., New York, 1969.


\bibitem{CP} C. De Concini and C. Procesi,
A characteristic free approach to invariant theory,  Adv. Math. 21
(1976), 330--354.

\bibitem{DKZ} M. Domokos, S. G.
Kuzmin, A. N. Zubkov, Rings of matrix invariants in positive
characteristic,  J. Pure Appl. Alg. 176 (2002), 61--80.

\bibitem{D1990} S. Donkin, The normality of closures of conjugacy
classes of matrices,  Inv. Math. 101 (1990), 717--736.

\bibitem{D1992} S. Donkin, Invariants of several
matrices,  Inv. Math.  110 (1992), 389--401.

\bibitem{m} H. Matsumura,
Commutative Ring Theory, Cambridge Univ. Press, 1986.

\bibitem{N} P. Newstead,  
Introduction to moduli problems and orbit spaces,
 Tata Inst. Lecture Notes, Springer-Verlag, 1978. 

\bibitem{Ri} D. R. Richman, The fundamental theorems of vector invariants, 
Adv. Math. 73 (1989), 43--78. 

\bibitem{R} M. Rosenlicht, Some basic theorems on algebraic groups,
Am. J. Math. 78 (1956), 401--443.

\bibitem{S} A. Seidenberg,
The hyperplane sections of normal varieties, Trans. Amer. Math.
Soc. 69 (1950), 357--386.

\bibitem{T} D. E. Taylor, 
The Geometry of the Classical
Groups, Heldermann Verlag, Berlin, 1992. 

\bibitem{W} H. Weyl, The Classical Groups --- Their Invariants and
Representations, Princeton University Press, 1946. 

\end{thebibliography}
\end{document}